\documentclass[12pt]{article}
\usepackage{amssymb,latexsym,amsmath}
\usepackage{graphicx,times}
\usepackage[margin=1.0in]{geometry}
\usepackage{multirow}
\newtheorem{The}{Theorem}[section]
\newtheorem{Def}{Definition}[section]

\newtheorem{cor}{Corollary}
\newtheorem{Ex}{Example}[section]
\numberwithin{equation}{section}
\begin{document}
\begin{center}
{\LARGE {\bf A Novel Third Order Numerical Method for Solving Volterra Integro-Differential Equations}}
\vskip 0.5cm
{\Large Sachin Bhalekar, Jayvant Patade\footnote{Corresponding author}}\\
\textit{Department of Mathematics, Shivaji University, Kolhapur - 416004, India.\\ Email:  sachin.math@yahoo.co.in, sbb\_maths@unishivaji.ac.in (Sachin Bhalekar), jayvantpatade1195@gmai.com (Jayvant Patade)}\\
\end{center}

\begin{abstract}
In this paper we introduce a numerical method for solving nonlinear Volterra  integro-differential equations. In the first step, we apply implicit trapezium rule to discretize the integral in given equation. Further, the  Daftardar-Gejji and Jafari technique (DJM) is used to find the unknown term on the right side. We derive existence-uniqueness theorem for such equations by using Lipschitz condition. We further present the  error,   convergence, stability and bifurcation analysis of the proposed method. We solve various types of equations using this method and compare the error with other numerical methods. It is observed that our method is more efficient than other numerical methods.
\end{abstract}
%\vskip 1cm
%\noindent

\section{Introduction}
The equations involving derivative of a dependent variable are frequently used to model the natural phenomena. The derivative of such function at a point $t_0$ can be approximated by using its values in a small neighborhood of $t_0$. Thus the derivative is a local operator which cannot model the memory and hereditary properties involved in the real world problems. To overcome this drawback, one has to include some nonlocal operator  in the given equation. Volterra \cite{VOL1, VOL2} suggested the use of an  integral operator (which is nonlocal) to model the problems in population dynamics. In his monograph \cite{VOL1}, he discussed the theory of integral, integro-differential and functional equations. Existence-uniqueness, stability and applications of integro-differential equations (IDE) is presented in a book by Lakshmikantham and Rao \cite{LAX}. Existence theory of nonlinear IDEs is also discussed in \cite{OREGAN}. Recently, Kostic \cite{Kostic}  discussed the theory of abstract Volterra integro-differential equations. 
\par A-stable linear multi-step methods to solve Volterra IDEs (VIDE) are proposed by Matthys in \cite{Matthys}. Brunner presented various numerical methods to solve VIDEs in  \cite{brunner1}. Day \cite{Day} used trapezoidal rule to devise a numerical method to solve nonlinear VIDEs. Linz \cite{Linz} derived fourth order numerical methods for such equations. Runge-Kutta method, predictor-corrector method and explicit multistep method are derived for IDEs by Wolfe and Phillips in \cite{Wolfe}. Dehghan and Salehi \cite{Dehghan} developed a numerical scheme based on the moving least square method for these equations. Singular IDEs in generalized Holder spaces are solved by using collocation and mechanical quadrature methods in \cite{Caraus}. Saeedi et al. \cite{Saeedi} described the operational Tau method for solving
 nonlinear VIDEs of the second kind.  The stability of numerical methods for VIDEs is studied in \cite{Brunner, Baker}.  A survey of numerical treatments of VIDEs is taken in \cite{Brunner3}.
\par The applications of integro-differential equations can be found in many fields such as mechanics and electromagnetic theory \cite{Bloom}, nuclear reactor \cite{Kastenberg, LAX}, visco-elasticity \cite{MacCamy, Narain}, heat conduction in materials with memory \cite{Hrusa, Pruss} and man-environment epidemics \cite{Capasso}.
\par In the present work we utilize a decomposition method proposed by Daftardar-Gejji and Jafari (DJM) \cite{GEJJI1} to generate a more accurate and faster numerical method for solving nonlinear VIDEs.The paper is organized as follows:
\par The preliminaries are given in section \ref{Pre}. A new numerical method is presented in Section \ref{nnm}. Analysis of this numerical method is given in Section \ref{analy}. Section \ref{example}  deals with different types of illustrative examples.  In Section \ref{Package}  we provide software package and conclusions are summarized in Section \ref{concl}.

\section{Preliminaries}\label{Pre}
\subsection{Basic Definitions and Results:}
In this section, we discuss some basic definitions and results\cite{Srivastava,LAX,Linz}.
\begin{Def}\label{1.11}
Let $y_j$  be the approximation to the exact value $y(x_j)$ obtained by a given method with step-size h. Then a method is said to be convergent if and only if
\begin{equation}
 \lim\limits_{h\rightarrow 0} \mid y(x_j)-y_j\mid\rightarrow 0,\quad j=1,2\cdots N.\label{1.1}
\end{equation}
\end{Def}
\begin{Def}\label{1.12}
A method is said to be of order p if p is the largest number for which there exists a finite constant $ \textsc{C}$ such that
\begin{equation}
  \mid y(x_j)-y_j\mid\leq \textsc{C}h^p,\quad j=1,2\cdots N.\label{1.2}
\end{equation}
\end{Def}

%\begin{Def}\label{1.5}
%Let $f\in C_\alpha $ and $\alpha \geq -1$, then the (left-sided) Riemann–Liouville integral of order $\mu, \mu> 0 $ is given by
%\begin{equation}
%  I^\mu f(t)=\frac{1}{\Gamma(\mu)} \int_{0}^t (t-\tau)^{\mu-1}f(\tau)d\tau,\quad t>0.
%\end{equation}
%\end{Def}

%\begin{Def}\label{1.6}
%The (left sided) Caputo fractional derivative of $f, f \in C_{-1}^m, m\in\mathbb{N}\cup\{0\}$, is defined as:
%\begin{eqnarray}
%D^\mu f(t)&=&\frac{d^m}{ dt^m} f(t),\quad \mu = m \nonumber\\
%&=& I^{m-\mu}\frac{d^m}{ dt^m} f(t),\quad {m-1} <\mu <m,\quad m\in \mathbb{N}.
%\end{eqnarray}
%\end{Def}

\begin{The}\label{1.3}
Consider the integral equation
\begin{equation}
 y(t) = f(t) + \int_0^x K(x,t, y(t))dt, \quad y(x_0)=y_0. \label{1.4}
\end{equation}
Assume that \\
(i) f is continuous in $0 \leq x \leq a$ ,\\
(ii) K(x, t, y) is continuous function for  $0 \leq t \leq x \leq a$ and $\parallel y\parallel < \infty$,\\
(iii) K(x, t, y) satisfies a Lipschitz condition
\begin{equation}
\parallel K(x, t, y_1)-K(x, t, y_2)\parallel \leq L \parallel y_1-y_2 \parallel
\end{equation}
for all  $0 \leq t \leq x \leq a$. Then the Eq. (\ref{1.4}) has a unique  solution.
\end{The}

\begin{The}\label{1.6}
Assume that $f\in C[I\times\mathbb{R}^n,\mathbb{R}^n]$, $K\in C[I\times I\times\mathbb{R}^n,\mathbb{R}^n]$ and $\int_{s}^x \mid K(t,s,y(s))\mid dt \leq N$, for $x_0\leq s \leq x \leq x_0 + a$, $y\in \Omega = \{\phi\in C[I,\mathbb{R}^n]: \phi(x_0) = x_0 \  and \  \mid\phi(x) - y_0\mid \leq b\}$. Then IVP (\ref{1.4}) possesses at least one solution.
\end{The}

\subsection{Existence and Uniqueness Theorem}\label{existnce}
In this subsection we prove existence and uniqueness theorem, which is generalization of Theorem \ref{1.3}.

\begin{The}\label{11}
Consider the Volterra integro-differential equation
\begin{eqnarray}
 y'(x) &=& f(x,y(x)) + \int_{x_0}^x K(x, t,y(t)) dt,\label{1}\\
\textrm{with initial condition}  \quad y(x_0) &=& y_0.\nonumber
\end{eqnarray}
Assume that  $f$ and $K$ are continuous and satisfy Lipschitz condition
\begin{eqnarray}
\parallel f(x, y_1)-f(x, y_2)\parallel &\leq& L_1 \parallel y_1-y_2 \parallel \\
\parallel K(x, t, y_1)-K(x, t, y_2)\parallel &\leq& L_2 \parallel y_1-y_2 \parallel
\end{eqnarray}
for every $\mid x-x_0\mid \leq a,  \mid t-x_0\mid \leq a, \parallel y_1 \parallel < \infty, \parallel y_2 \parallel < \infty\quad \textrm{and}\quad a>0$.
Then the Eq. (\ref{1}) has unique solution.
\end{The}
\textbf{Proof:}
Integrating Eq.(\ref{1}) and using $ y(x_0) = y_0$, we get
\begin{eqnarray}
 y(x) & = & y_0 + \int_{x_0}^x f(x, y(x)) dx + \int_{x_0}^x \left(\int_{x_0}^x K(x, t, y(t))dt\right)dx\nonumber\\
\Rightarrow  y(x) & = & y_0 + \int_{x_0}^x \left(f(x, y(x)) + \int_{x_0}^x K(x, t, y(t))dt\right)dx\nonumber\\
\Rightarrow y(x) &=& y_0 + \int_{x_0}^x G(x, y(x)) dx,\label{9.1}\\
\textrm{where}\quad G(x,y(x)) &=& f(x, y(x)) + \int_{x_0}^x K(x, t, y(t))dt.\nonumber
\end{eqnarray}
The Eq. (\ref{9.1}) is of the form (\ref{1.4}).\\
 We have following observations\\
(i) $y_0$  is continuous  because it is constant,\\
(ii) Kernel $G$ is continuous for $0 \leq x \leq a$, because $f$ and $K$ are continuous in the same domain, \\
\begin{eqnarray*}
(iii) \parallel G(x, y_1)-G(x, y_2)\parallel &=&  \parallel f(x, y_1(x)) + \int_{x_0}^x K(x, t, y_1(t))dt - f(x, y_2(x))\\
&& - \int_{x_0}^x K(x, t, y_2(t))dt\parallel \\
&\leq&\parallel f(x, y_1(x))-f(x, y_2(x))\parallel + \parallel \int_{x_0}^x K(x, t, y_1(t))dt\\
&& - \int_{x_0}^x K(x, t, y_2(t))dt\parallel \\
&\leq& L_1 \parallel y_1-y_2 \parallel +  a L_2 \parallel y_1-y_2 \parallel \quad (\because \mid x-x_0\mid \leq a)\\
&\leq& (L_1 + aL_2)\parallel y_1-y_2 \parallel
\end{eqnarray*}
$\therefore G$ satisfy a Lipschitz condition.\\
$\Rightarrow$ All the conditions of the Theorem \ref{1.3} are satisfies.\\
Hence the  Eq.(\ref{1}) has unique solution.

\subsection{Daftardar-Gejji and Jafari Method}\label{djm}
A new iterative method was introduced by Daftardar-Gejji and Jafari (DJM) \cite{GEJJI1} in 2006 for solving nonlinear functional equations. The DJM has been used to solve a variety of equations such as  fractional differential equations \cite{GEJJI3}, partial differential equations \cite{GEJJI4}, boundary value problems \cite{GEJJI5, Mohyud}, evolution equations \cite{GEJJI6} and system of nonlinear functional equations \cite{GEJJI8}. The method is successfully employed to solve Newell-Whitehead-Segel equation \cite{pb1}, fractional-order logistic equation \cite{GEJJI7} and some nonlinear  dynamical systems \cite{GEJJI9} also. Recently DJM has been used to generate new numerical methods \cite{GEJJI10,GEJJI11,pb2} for solving differential equations.
In this section we describe DJM which is very useful for solving the equations of the form
\begin{equation}
	u= g + L(u) + N(u),\label{2.1}
\end{equation}
 where $g$ is a given function, $L$ and $N$ are linear and nonlinear operators respectively. DJM  provides  the  solution to Eq.(\ref{2.1}) in the form of series
 \begin{equation}
	u= \sum_{i=0}^\infty u_i. \label{2.2}
\end{equation}
where
\begin{eqnarray}
u_0 &= & g, \nonumber\\
u_{m+1} &=& L(u_m) + G_m,\quad m=0,1, 2, \cdots,\label{2.6}\\
\textrm{and}\quad G_m &=& N\left(\sum_{j=0}^m u_m\right)- N\left(\sum_{j=0}^{m-1} u_m\right), m\geq 1.\nonumber
 \end{eqnarray}
The $k$-term approximate solution is given by
\begin{eqnarray}
u = \sum_{i=0}^{k-1} u_i
 \end{eqnarray}
for suitable integer $k$.

 The following convergence results for DJM are described in \cite{CONV}.
\begin{The}
If $N$ is $C^{(\infty)}$ in a neighborhood of $u_0$ and $\left\|N^{(n)}(u_0) \right\| \leq L$, for any $n$ and for some real $L>0$ and $\left\|u_i\right\| \leq M <\frac{1}{e}$, $i=1,2,\cdots,$ then the series $\sum_{n=0}^{\infty} G_n$ is absolutely convergent to $N$ and moreover,
\begin{equation*}
\left\|G_n\right\| \leq L M^n e^{n-1} (e-1), \quad n=1,2,\cdots.
\end{equation*}
\end{The}

\begin{The}
If $N$ is $C^{(\infty)}$ and $\left\|N^{(n)}(u_0) \right\| \leq M \leq e^{-1}$, $\forall n$, then the series $\sum_{n=0}^{\infty} G_n$ is absolutely convergent to $N$.
\end{The}

\section{Numerical Method}\label{nnm}
In this section we present a numerical method based on DJM to solve Volterra integro-differential equation. 

Integrating Eq. (\ref{1}) from $x=x_j $ to $x=x_j+h$, we get

\begin{equation}
y({x_j+h}) = y(x_j) + \int_{x_j}^{x_j+h} f(x,y(x)) dx + \int_{x_j}^{x_j+h}  \int_{x_0}^x K(x,t, y(t))dt dx \label{2}
\end{equation}
Applying trapezium formula \cite{Jain} to evaluate integrals on right of Eq. (\ref{2}), we get
\begin{eqnarray}
y({x_j+h}) &=& y(x_j) + \frac{h}{2} f(x_j, y_j) + \frac{h^2}{4} \left( K(x_j, x_0, y_0) + K(x_j, x_j, y_j) + K(x_{j+1}, x_0, y_0)\right)\nonumber\\
&& + \frac{h^2}{2}\left( \sum_{i=1}^{j-1} K(x_j, x_i, y_i) +  \sum_{i=1}^j K(x_{j+1}, x_i, y_i)\right) + \frac{h}{2} f(x_{j+1}, y_{j+1}) \nonumber\\
&&+ \frac{h^2}{4}  K(x_{j+1}, x_{j+1}, y_{j+1}) + O(h^3).\label{3}
\end{eqnarray}
If $y_j$ is an approximation to $y(x_j)$, then approximate solution is given by
\begin{eqnarray}
y_{j+1} &=& y_j + \frac{h}{2} f(x_j, y_j) + \frac{h^2}{4} \left( K(x_j, x_0, y_0) + K(x_j, x_j, y_j) + K(x_{j+1}, x_0, y_0)\right)\nonumber\\
&& + \frac{h^2}{2}\left( \sum_{i=1}^{j-1} K(x_j, x_i, y_i) +  \sum_{i=1}^j K(x_{j+1}, x_i, y_i)\right) + \frac{h}{2} f(x_{j+1}, y_{j+1}) \nonumber\\
&& + \frac{h^2}{4}  K(x_{j+1}, x_{j+1}, y_{j+1}).\label{4.1}
\end{eqnarray}
Equation (\ref{4.1}) is of the form (\ref{2.1}),
where
\begin{eqnarray*}
u &=& y_{j+1},\nonumber\\
g &=& y_j + \frac{h}{2} f(x_j, y_j) + \frac{h^2}{4} \left( K(x_j, x_0, y_0) + K(x_j, x_j, y_j) + K(x_{j+1}, x_0, y_0)\right)\\
&& + \frac{h^2}{2}\left( \sum_{i=1}^{j-1} K(x_j, x_i, y_i) +  \sum_{i=1}^j K(x_{j+1}, x_i, y_i)\right)\\
N(u) &=& \frac{h}{2} f(x_{j+1}, y_{j+1}) + \frac{h^2}{4}  K(x_{j+1}, x_{j+1}, y_{j+1}).
\end{eqnarray*}
Applying DJM to equation (\ref{4.1}), we obtain 3-term solution as
\begin{eqnarray*}
u &=& u_0 + u_1 + u_2\\
&=& u_0 + N(u_0) + N(u_0+u_1) - N(u_0)\\
&=& u_0 + N(u_0+u_1)\\
&=& u_0 + N(u_0+N(u_0)).
\end{eqnarray*}
That is
\begin{eqnarray}
y_{j+1} &=& y_j + \frac{h}{2} f(x_j, y_j) + \frac{h^2}{4} \left( K(x_j, x_0, y_0) + K(x_j, x_j, y_j) + K(x_{j+1}, x_0, y_0)\right)\nonumber\\
&& + \frac{h^2}{2}\left( \sum_{i=1}^{j-1} K(x_j, x_i, y_i) +  \sum_{i=1}^j K(x_{j+1}, x_i, y_i)\right)\nonumber\\
&& + N\left(y_j + \frac{h}{2} f(x_j, y_j) + \frac{h^2}{4} \left( K(x_j, x_0, y_0) + K(x_j, x_j, y_j) + K(x_{j+1}, x_0, y_0)\right)\right.\nonumber\\
&& \left. + \frac{h^2}{2}\left( \sum_{i=1}^{j-1} K(x_j, x_i, y_i) +  \sum_{i=1}^j K(x_{j+1}, x_i, y_i)\right)\right) + N\left(y_j + \frac{h}{2} f(x_j, y_j)\right.\nonumber\\
&&\left. + \frac{h^2}{4} \left( K(x_j, x_0, y_0) + K(x_j, x_j, y_j) + K(x_{j+1}, x_0, y_0)\right)\right.\nonumber\\
&&\left. + \frac{h^2}{2}\left( \sum_{i=1}^{j-1} K(x_j, x_i, y_i) +  \sum_{i=1}^j K(x_{j+1}, x_i, y_i)\right)\right).
\end{eqnarray}
or
\begin{eqnarray}
y_{j+1} &=& y_j + \frac{h}{2} f(x_j, y_j) + \frac{h^2}{4} \left( K(x_j, x_0, y_0) + K(x_j, x_j, y_j) + K(x_{j+1}, x_0, y_0)\right)\nonumber\\
&& + \frac{h^2}{2}\left( \sum_{i=1}^{j-1} K(x_j, x_i, y_i) +  \sum_{i=1}^j K(x_{j+1}, x_i, y_i)\right)\nonumber\\
&& + \frac{h}{2}f\left(x_{j+1}, y_j + \frac{h}{2} f(x_j, y_j) + \frac{h^2}{4} \left( K(x_j, x_0, y_0) + K(x_j, x_j, y_j) + K(x_{j+1}, x_0, y_0)\right)\right.\nonumber\\
&& \left. + \frac{h^2}{2}\left( \sum_{i=1}^{j-1} K(x_j, x_i, y_i) +  \sum_{i=1}^j K(x_{j+1}, x_i, y_i)\right)+ \frac{h}{2}f\left(x_{j+1}, y_j + \frac{h}{2} f(x_j, y_j) \right.\right.\nonumber\\
&&\left.\left. + \frac{h^2}{4} \left( K(x_j, x_0, y_0) + K(x_j, x_j, y_j) + K(x_{j+1}, x_0, y_0)\right)\right.\right.\nonumber\\
&&\left.\left. + \frac{h^2}{2}\left( \sum_{i=1}^{j-1} K(x_j, x_i, y_i) +  \sum_{i=1}^j K(x_{j+1}, x_i, y_i)\right)\right)\right.\nonumber\\
&&\left. + \frac{h^2}{4} K\left(x_{j+1},x_{j+1},y_j + \frac{h}{2} f(x_j, y_j) + \frac{h^2}{4} \left( K(x_j, x_0, y_0) + K(x_j, x_j, y_j) + K(x_{j+1}, x_0, y_0)\right)\right.\right.\nonumber\\
&& \left.\left.+ \frac{h^2}{2}\left( \sum_{i=1}^{j-1} K(x_j, x_i, y_i) +  \sum_{i=1}^j K(x_{j+1}, x_i, y_i)\right)\right)\right)\nonumber\\
&& + \frac{h^2}{4} K\left(x_{j+1},x_{j+1},y_j + \frac{h}{2} f(x_j, y_j) + \frac{h^2}{4} \left( K(x_j, x_0, y_0) + K(x_j, x_j, y_j) + K(x_{j+1}, x_0, y_0)\right)\right.\nonumber\\
&& \left. + \frac{h^2}{2}\left( \sum_{i=1}^{j-1} K(x_j, x_i, y_i) +  \sum_{i=1}^j K(x_{j+1}, x_i, y_i)\right)+ \frac{h}{2}f\left(x_{j+1}, y_j + \frac{h}{2} f(x_j, y_j) \right.\right.\nonumber\\
&&\left.\left. + \frac{h^2}{4} \left( K(x_j, x_0, y_0) + K(x_j, x_j, y_j) + K(x_{j+1}, x_0, y_0)\right)\right.\right.\nonumber\\
&&\left.\left. + \frac{h^2}{2}\left( \sum_{i=1}^{j-1} K(x_j, x_i, y_i) +  \sum_{i=1}^j K(x_{j+1}, x_i, y_i)\right)\right)\right.\nonumber\\
&&\left. + \frac{h^2}{4} K\left(x_{j+1},x_{j+1},y_j + \frac{h}{2} f(x_j, y_j) + \frac{h^2}{4} \left( K(x_j, x_0, y_0) + K(x_j, x_j, y_j) + K(x_{j+1}, x_0, y_0)\right)\right.\right.\nonumber\\
&& \left.\left.+ \frac{h^2}{2}\left( \sum_{i=1}^{j-1} K(x_j, x_i, y_i) +  \sum_{i=1}^j K(x_{j+1}, x_i, y_i)\right)\right)\right)\label{5}
\end{eqnarray}
If we set
\begin{eqnarray}
M_1 &=& y_j + \frac{h}{2} f(x_j, y_j) + \frac{h^2}{4} \left( K(x_j, x_0, y_0) + K(x_j, x_j, y_j) + K(x_{j+1}, x_0, y_0)\right)\nonumber\\
&& + \frac{h^2}{2}\left( \sum_{i=1}^{j-1} K(x_j, x_i, y_i) +  \sum_{i=1}^j K(x_{j+1}, x_i, y_i)\right), \label{6}\\
M_2 &=& M_1 +  \frac{h}{2} f(x_{j+1}, M_1) +  \frac{h^2}{4}  K(x_{j+1}, x_{j+1}, M_1)\label{7}
\end{eqnarray}
then equation (\ref{5}) becomes
\begin{eqnarray}
y_{j+1}=M_1 + \frac{h}{2} f(x_{j+1}, M_2) +  \frac{h^2}{4}  K(x_{j+1}, x_{j+1}, M_2).\label{9}
\end{eqnarray}

\section{Analysis of Numerical Method}\label{analy}
\subsection{Error Analysis }
\begin{The}\label{13}
The numerical method (\ref{9}) is of third order.
\end{The}
\textbf{Proof:}
Suppose  $y_{j+1}$ is an approximation to $y(x_{j+1})$. By using Eq.(\ref{3}) and (\ref{9}), we obtain
\begin{eqnarray*}
\mid y(x_{j+1})-y_{j+1}\mid &=& \mid\frac{h}{2} f(x_{j+1}, y_{j+1}) + \frac{h^2}{4}  K(x_{j+1}, x_{j+1}, y_{j+1}) + O(h^3)\\
&&-\frac{h}{2} f(x_{j+1}, M_2)-\frac{h^2}{4}  K(x_{j+1}, x_{j+1}, M_2)\mid\\
&\leq& \frac{h}{2}\mid f(x_{j+1}, y_{j+1})- f(x_{j+1}, M_2)\mid + \frac{h^2}{4}\mid  K(x_{j+1}, x_{j+1}, y_{j+1})\\
&&- K(x_{j+1}, x_{j+1}, M_2)\mid+ O(h^3)\\
&\leq& \left(\frac{h}{2}L_1 + \frac{h^2}{4}L_2\right)\mid y_{j+1}-M_2\mid + O(h^3).
\end{eqnarray*}
Using Eq.(\ref{7}) and (\ref{9}), we get
\begin{eqnarray*}
 \mid y(x_{j+1})-y_{j+1}\mid &\leq& \left(\frac{h}{2}L_1 + \frac{h^2}{4}L_2\right)\mid M_1 + \frac{h}{2} f(x_{j+1}, M_2) +  \frac{h^2}{4}  K(x_{j+1}, x_{j+1}, M_2)\\
&&- M_1 - \frac{h}{2} f(x_{j+1}, M_1)- \frac{h^2}{4}  K(x_{j+1}, x_{j+1}, M_1)\mid+ O(h^3)\\
&\leq& \left(\frac{h}{2}L_1 + \frac{h^2}{4}L_2\right)\left(\mid \frac{h}{2} f(x_{j+1}, M_2)-\frac{h}{2} f(x_{j+1}, M_1)\mid\right.\\
&&\left. + \mid \frac{h^2}{4}  K(x_{j+1}, x_{j+1}, M_2)-\frac{h^2}{4}  K(x_{j+1}, x_{j+1}, M_1)\mid\right)+ O(h^3).
\end{eqnarray*}
Using Eq.(\ref{6}) and (\ref{7}) ,we get
\begin{eqnarray*}
 \mid y(x_{j+1})-y_{j+1}\mid &\leq& \left(\frac{h}{2}L_1 + \frac{h^2}{4}L_2\right)^2\mid \frac{h}{2} f(x_{j+1}, M_1) +  \frac{h^2}{4}  K(x_{j+1}, x_{j+1}, M_2)+ O(h^3)\\
&\leq& h^3 \left(\left(\frac{1}{2}L_1 + \frac{h}{4}L_2\right)^2\mid \frac{1}{2} f(x_{j+1}, M_1) +\frac{h}{4}  K(x_{j+1}, x_{j+1}, M_2)\mid\right) + O(h^2).
\end{eqnarray*}
$\Rightarrow$ The numerical method (\ref{9}) is of third order.

\begin{cor}\label{14}
The numerical method (\ref{9}) is convergent.
\end{cor}
\textbf{Proof:} By Theorem \ref{13} and definition (\ref{1.11}) the  numerical method (\ref{9}) is convergent.

\subsection{Stability Analysis of Numerical Method}\label{sability}
Consider the test equation 
\begin{equation}
 y' = \alpha y + \beta\int_{x_0}^x y(t)dt. \label{15}
\end{equation}
 Applying numerical method (\ref{9}) to this equation, we get
 \begin{equation}
y_{j+1}= M_1 + \frac{h\alpha}{2}M_2 + \frac{h^2\beta}{4}M_2,\label{18}
 \end{equation}
where 
\begin{eqnarray}
M_1 &=& \frac {h^2\beta}{2}y_0 + h^2\beta \sum_{i=1}^{j-1} y_i +\left(1+ \frac {h\alpha}{2} + \frac {3h^2\beta}{4}\right)y_j ,\\
M_2 &=& \left(1+ \frac {h\alpha}{2} + \frac {h^2\beta}{4}\right)M_1.
\end{eqnarray}
 If we set $u=h\alpha$ and $v=h^2\beta$ then Eq. (\ref{18}) can be written as
\begin{eqnarray}
y_{j+1} &=& \left(\frac{v}{2}+\frac{u v}{4}+\frac{u^2 v}{8}+\frac{v^2}{8}+\frac{u v^2}{8}+\frac{v^3}{32}\right)y_0\nonumber\\
&&+ \left(v+\frac{u v}{2}+\frac{u^2 v}{4}+\frac{v^2}{4}+\frac{u v^2}{4}+\frac{v^3}{16}\right)\sum_{i=1}^{j-1} y_i \nonumber\\
&&+ \left(1+u+\frac{u^2}{2}+\frac{u^3}{8}+v+\frac{3 u v}{4}+\frac{5 u^2 v}{16}+\frac{v^2}{4}+\frac{7 u v^2}{32}+\frac{3 v^3}{64}\right)y_j.\label{16}
\end{eqnarray}
Simplifying Eq.(\ref{16}), we get
\begin{equation}
y_j = b_1 y_{j-1} + b_2 y_{j-2},\label{17}
 \end{equation}
where
\begin{eqnarray}
b_1 &=& 2+u+\frac{u^2}{2}+\frac{u^3}{8}+v+\frac{3 u v}{4}+\frac{5 u^2 v}{16}+\frac{v^2}{4}+\frac{7 u v^2}{32}+\frac{3 v^3}{64},\\
\textrm{and}\quad b_2 &=& -1-u-\frac{u^2}{2}-\frac{u^3}{8}-\frac{u v}{4}-\frac{u^2 v}{16}+\frac{u v^2}{32}+\frac{v^3}{64}.
\end{eqnarray}
The characteristic equation \cite{Elaydi} of the difference  Eq. (\ref{17}) is
\begin{equation}
r^2 - b_1 r -b_2 = 0.
\end{equation}
The characteristic roots are given by
\begin{equation}
 r_1 =\frac{b_1 + \sqrt{b_1^2+4 b_2}}{2}, \quad r_2 = \frac{b_1 - \sqrt{b_1^2+4 b_2}}{2}. 
\end{equation}
The zero solution of system (\ref{17}) is asymptotically stable if $\mid r_1 \mid< 1 $ and  $\mid r_2 \mid< 1 $.
Thus the stability region of numerical method (\ref{9}) is given by  $\mid r_1 \mid< 1 $ and  $\mid r_2 \mid< 1 $ as shown in  Fig.1.\\
\\
\begin{center}
\begin{tabular}{c}
\includegraphics[scale=0.8]{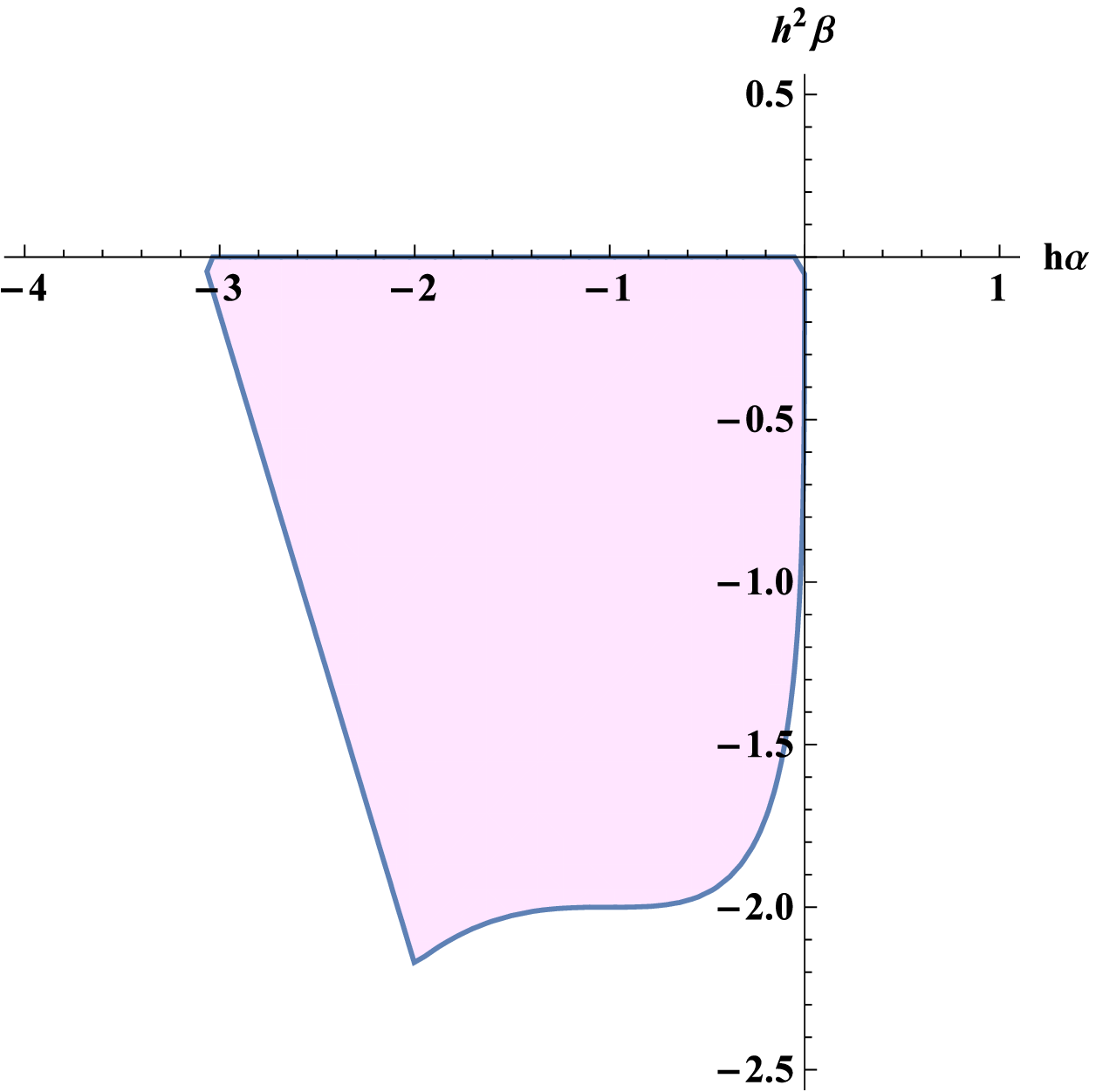}\\
Fig.1 Stability region for  numerical method (\ref{9})\\
\end{tabular}
\end{center}

\subsection{Bifurcations Analysis}\label{bianaly}
We discuss bifurcation analysis for the following test equation given in \cite{Edwards}
\begin{equation}
 y' = -\int_{x_0}^x e^{-\alpha(x-t)} y(t)dt. \label{25}
\end{equation}
As discussed in \cite{Edwards}, Eq.(\ref{25}) with convolution kernel and with exponentially fading memory is practically important. The qualitative behavior of the Eq.(\ref{25}) is described as below \cite{Edwards}:\\
\textbf{Table 1: Bifurcation in IDE (\ref{25}).}\\
\begin{tabular}{|c|c|}
\hline
\textbf{Region} & \textbf{Qualitative behavior} \\
\hline
$\alpha \geq 2$  & Solutions are converging to zero without oscillation  \\
\hline
$0 <\alpha < 2$ & Solutions are damped oscillations converging to zero \\
\hline
$-2 <\alpha < 0$ & Divergent solutions exhibiting  unbounded oscillations \\
\hline 
 $\alpha \leq -2 $ & Solutions become unbounded without oscillations \\
\hline 
\end{tabular}\\

 Now we apply the numerical scheme (\ref{9}) to  Eq.(\ref{25}) and discuss the bifurcations.
 \begin{equation}
 y_{j+1} = -a_0 e^{-\alpha h j}y_0 - a_1 \sum_{i=1}^{j-1}  e^{-\alpha h (j-1)}y_i + a_2 y_j,
 \end{equation}
where
\begin{eqnarray*}
a_0 &=& \frac{h^2}{4} \left(1-\frac{h^2}{4} + \frac{h^4}{16}\right)\left(1+e^{-\alpha h}\right),\\
a_1 &=& \frac{h^2}{2} \left(1-\frac{h^2}{4} + \frac{h^4}{16}\right)\left(1+e^{-\alpha h}\right),\\
a_2 &=& \left(1-\frac{h^2}{4} + \frac{h^4}{16}\right) \left(1-\frac{h^2}{2} + \frac{h^2}{4} e^{-\alpha h}\right).
\end{eqnarray*}
Eq. (\ref{25}) can be written as
\begin{equation}
y_{j+2} - b_1 y_{j+1} + b_2 y_j = 0,  \label{26}
\end{equation}
where
\begin{eqnarray*}
b_1 &=& e^{-\alpha h}+a_2,\\
b_2 &=& (a_1 + a_2)e^{-\alpha h}.
\end{eqnarray*} 
The characteristic equation of the difference Eq. (\ref{26}) is
\begin{equation}
r^2-b_1 r + b_2 = 0.
\end{equation}
The roots of this equation are
\begin{equation}
r = \frac{b_1 \pm \sqrt{b_1^2 - 4b_2}}{2}.\label{27}
\end{equation}

\begin{The}\cite{Elaydi}\label{29}
All solutions of the difference equation
\begin{equation}
y_{j+2} + p_1 y_{j+1} + p_2 y_j = 0,\label{28}
\end{equation}
converge to zero (i.e.,the zero solution is asymptotically stable) if and only if $max\{|r_1|, |r_2|\} < 1$, where $r_1,r_2$ are the roots of the characteristic equation $r^2+p_1r+p_2=0$ of Eq. (\ref{28}).
\end{The}
Define
\begin{eqnarray*}
\alpha_0 &=& \frac{1}{h}\text{Log}\left(1+\frac{h^6}{64}\right),\\
\alpha_1 &=& \frac{1}{h}\text{Log}\left(\frac{2 \left(128+160 h^2-32 h^4+8 h^6-h^8+8 h \sqrt{1024+8 h^6-2 h^8}\right)}{\left(-4+h^2\right)^2 \left(16-4h^2+h^4\right)}\right)\quad \textrm{and}\\
\alpha_2 &=& \frac{1}{h}\text{Log}\left(-\frac{2 \left(-128-160 h^2+32 h^4-8 h^6+h^8+8 h \sqrt{1024+8 h^6-2 h^8}\right)}{\left(-4+h^2\right)^2\left(16-4 h^2+h^4\right)}\right).
\end{eqnarray*}
Using Theorem \ref{29} we have the following observations:\\
\textbf{Case 1:} The characteristics roots of the Eq. (\ref{27}) are real and distinct
\begin{eqnarray*}
r_1 &=& \frac{1}{2} \left(1+e^{-h \alpha }-\frac{h^2}{2}-\frac{1}{2} e^{-h \alpha } h^2+\frac{h^4}{8}+\frac{1}{8} e^{-h \alpha } h^4-\frac{h^6}{64}-\frac{1}{32}e^{-h \alpha } h^6 + \right.\\
&&\left.\sqrt{-4 e^{-h \alpha } \left(1+\frac{h^6}{64}\right)+\left(1+e^{-h \alpha }-\frac{h^2}{2}-\frac{1}{2} e^{-h \alpha } h^2+\frac{h^4}{8}+\frac{1}{8}e^{-h \alpha } h^4-\frac{h^6}{64}-\frac{1}{32} e^{-h \alpha } h^6\right)^2}\right)\\
r_2 &=& \frac{1}{2} \left(1+e^{-h \alpha }-\frac{h^2}{2}-\frac{1}{2} e^{-h \alpha } h^2+\frac{h^4}{8}+\frac{1}{8} e^{-h \alpha } h^4-\frac{h^6}{64}-\frac{1}{32}e^{-h \alpha } h^6 - \right.\\
&&\left.\sqrt{-4 e^{-h \alpha } \left(1+\frac{h^6}{64}\right)+\left(1+e^{-h \alpha }-\frac{h^2}{2}-\frac{1}{2} e^{-h \alpha } h^2+\frac{h^4}{8}+\frac{1}{8}e^{-h \alpha } h^4-\frac{h^6}{64}-\frac{1}{32} e^{-h \alpha } h^6\right)^2}\right),
\end{eqnarray*}
when\\
\textbf{Case 1.1:} $\alpha > \alpha_1$. \\
In this case the characteristic roots are of magnitude less than unity (cf. Region I in Fig.2). Therefore, all the solutions of difference  Eq. (\ref{26}) converge to zero without oscillations.\\
\textbf{Case 1.2:} $\alpha <\alpha_2$.\\
In this case the characteristic roots are of magnitude greater than unity. This is shown in Fig.2 as Region IV. In this region, all the solutions will tend to infinity without oscillations.\\
\textbf{Case 2:} The characteristics roots are real and equal
\begin{equation*}
r=r_1=r_2=\frac{1}{2} \left(1+e^{-h \alpha }-\frac{h^2}{2}-\frac{1}{2} e^{-h \alpha } h^2+\frac{h^4}{8}+\frac{1}{8} e^{-h \alpha} h^4-\frac{h^6}{64}-\frac{1}{32} e^{-h \alpha } h^6\right),
\end{equation*}
when\\
\textbf{Case 2.1:} $\alpha = \alpha_1$.\\
In this case the characteristic root $r\in(0,1)$. Hence the solutions of Eq. (\ref{26}) tend to zero without oscillations as in case 1.1. This case is shown in Fig.2 as boundary line of Region I and  Region II.\\
\textbf{Case 2.2:} $\alpha = \alpha_2$. \\
In this case the characteristic root $r>1$ and the solutions are unbounded (cf. boundary line of Region III and  Region IV in Fig.2).\\
\textbf{Case 3:} The characteristics roots are complex
\begin{eqnarray*}
r_1 &=& \frac{1}{2} \left(1+e^{-h \alpha }-\frac{h^2}{2}-\frac{1}{2} e^{-h \alpha } h^2+\frac{h^4}{8}+\frac{1}{8} e^{-h \alpha } h^4-\frac{h^6}{64}-\frac{1}{32}
e^{-h \alpha } h^6\right.\\
&&\left.+ i \sqrt{4 e^{-h \alpha } \left(1+\frac{h^6}{64}\right)-\left(1+e^{-h \alpha }-\frac{h^2}{2}-\frac{1}{2} e^{-h \alpha } h^2+\frac{h^4}{8}+\frac{1}{8}
e^{-h \alpha } h^4-\frac{h^6}{64}-\frac{1}{32} e^{-h \alpha } h^6\right)^2}\right)\\
r_2 &=&  \frac{1}{2} \left(1+e^{-h \alpha }-\frac{h^2}{2}-\frac{1}{2} e^{-h \alpha } h^2+\frac{h^4}{8}+\frac{1}{8} e^{-h \alpha } h^4-\frac{h^6}{64}-\frac{1}{32}
e^{-h \alpha } h^6\right.\\
&&\left.- i \sqrt{4 e^{-h \alpha } \left(1+\frac{h^6}{64}\right)-\left(1+e^{-h \alpha }-\frac{h^2}{2}-\frac{1}{2} e^{-h \alpha } h^2+\frac{h^4}{8}+\frac{1}{8}
e^{-h \alpha } h^4-\frac{h^6}{64}-\frac{1}{32} e^{-h \alpha } h^6\right)^2}\right),
\end{eqnarray*}
when $\alpha_2 < \alpha <\alpha_1$.\\
Note that $\mid r_1\mid= \mid r_2\mid$.\\
\textbf{Case 3.1:} If $\alpha_0 < \alpha <\alpha_1$ then $\mid r_1\mid<1$. The solutions of Eq. (\ref{26}) are oscillatory and converge to zero (Region II in Fig.2). \\
\textbf{Case 3.2:} If $\alpha=\alpha_0$ then  $\mid r_1\mid=1$ and bifurcation occurs. This is shown by the boundary line of Region II and  Region III in Fig.2.\\
\textbf{Case 3.3:} If $\alpha_2 < \alpha < \alpha_0$\\
then $\mid r_1\mid > 1$. This implies that the solutions of Eq. (\ref{26}) are oscillatory and diverge to infinity (Region III in Fig.2).
 
 \begin{center}
 \begin{tabular}{c}
 \includegraphics[scale=1.0]{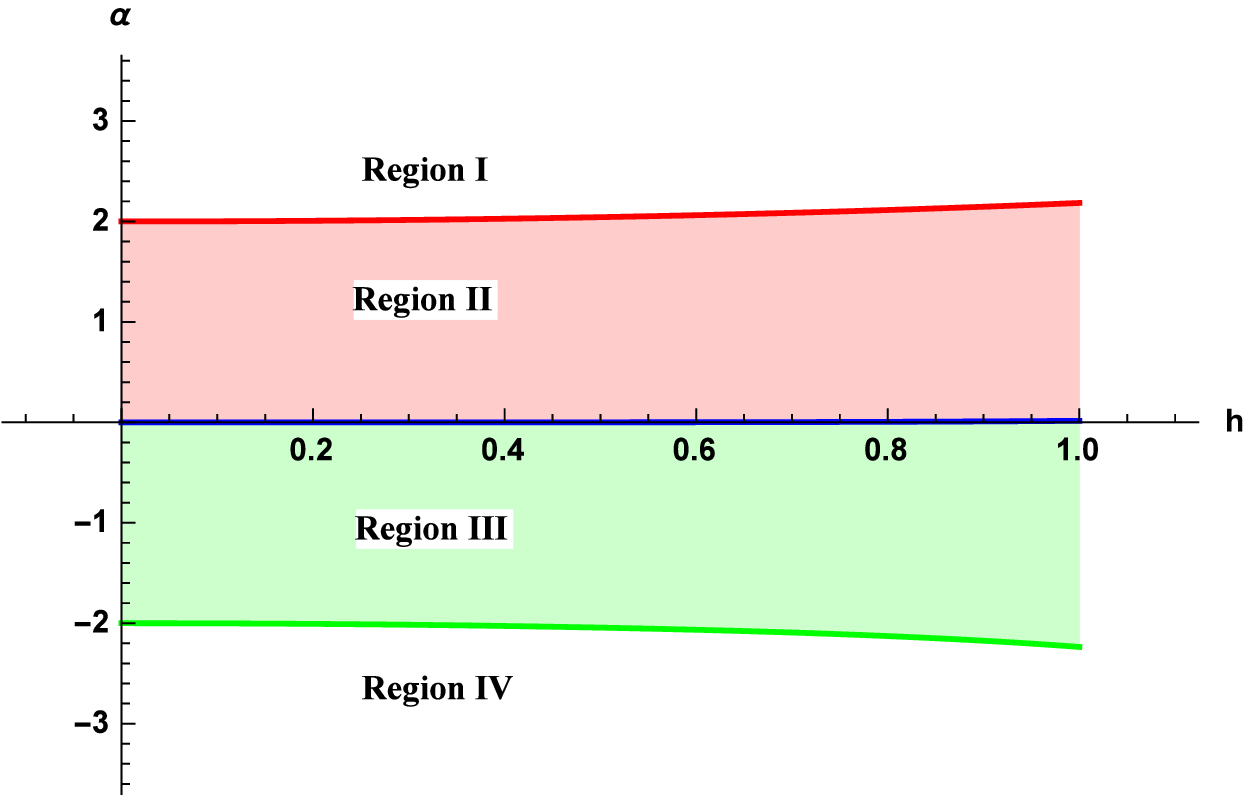}\\
 Fig.2: Bifurcation digram for Eq. (\ref{26}).\\
 \end{tabular}
 \end{center}

Note: It is clear from Table 1 and Fig.2 that the bifurcations in numerical scheme  (\ref{26}) coincide with those in IDE (\ref{25}) for $h=0$.

\section{Illustrative examples}\label{example}
\begin{Ex} \cite{Day,Wolfe}\label{ex1}
Consider the Volterra integro-differential equation
\begin{equation}
y'(x) = 1 + 2x - y(x) + \int_{x_0}^x x(1+2x)e^{t(x-t)}y(t)dt,\quad y(0) = 1.\label{19}
\end{equation}
\end{Ex}
The exact solution of Eq.(\ref{19}) is $y(x) = e^{x^2}$. We compare errors obtained using different numerical methods given in \cite{Day,Wolfe}, viz. third order Runge-Kutta method, predictor-corrector method, multi-step method and Day method with our method for $h=0.1$ and  $h=0.025$  in Table 2 and Table 3 respectively. It is observed that the error in our method is very less as compared with other numerical methods.\\
Note that the error values in Table 2 and Table 3 for all other methods are taken from \cite{Day,Wolfe}.
\newpage
\textbf{Table 2: Comparison of numerical results of example 6.1 with}  $\textbf{h=0.1}$ \\
\begin{tabular}{|c|c|c|c|c|c|c|}
\hline
\textbf{x} & \textbf{Our method} & \textbf{Runge-Kutta}  & \textbf{Predictor-} & \textbf{Multi-step} & \textbf{Day method}  \\
\multirow{2}{*}
 &$O(h^3)$  &  $O(h^3)$ & \textbf{Corrector} $O(h^3)$ & $O(h^3)$&   \\
\hline
  0.1 & $2.3\times10^{-5}$ & $ 5.5\times10^{-4}$ & $ 5.0\times10^{-5}$ & $5.0\times10^{-5}$& $5.4\times10^{-4}$ \\
\hline
  0.2 & $1.2 \times10^{-4}$ & $ 1.2\times10^{-3}$ & $ 5.2\times10^{-5}$& $8.1\times10^{-4}$& $1.1\times10^{-3}$ \\
\hline
  0.3 & $3.1\times10^{-4}$ & $ 1.9\times10^{-3}$  & $2.7\times10^{-4}$ & $7.5\times10^{-4}$& $1.9\times10^{-3}$ \\
\hline
  0.4 & $6.2\times10^{-3}$ & $2.8\times10^{-3}$  & $6.6\times10^{-4}$ & $2.1\times10^{-3}$ & $2.8\times10^{-3}$ \\
\hline
  0.5 & $1.0\times10^{-3}$ & $3.9\times10^{-3}$ & $1.2\times10^{-3}$ & $2.5\times10^{-3}$ & $4.0\times10^{-3}$ \\
\hline
 0.6 & $1.7\times10^{-3}$ & $ 5.8\times10^{-3}$  & $ 2.0\times10^{-3}$ & $4.7\times10^{-3}$& $5.6\times10^{-3}$ \\
\hline
  0.7 & $2.8\times10^{-3}$ & $7.3 \times10^{-3} $ & $3.1\times10^{-3}$ & $6.0\times10^{-3}$& $7.9\times10^{-3}$ \\
\hline
  0.8 & $4.4\times10^{-3}$ & $1.0 \times10^{-2} $ & $4.8\times10^{-3}$ & $1.0\times10^{-2}$& $1.1\times10^{-2}$ \\
\hline
  0.9 & $6.8\times10^{-3}$ & $1.4 \times10^{-2} $ & $7.3\times10^{-3}$ & $1.4\times10^{-2}$& $1.5\times10^{-2}$\\
\hline
  1 & $ 1.0\times10^{-2}$ &  $2.0\times10^{-2} $ & $1.1\times10^{-2}$ & $2.1\times10^{-2}$& $2.2\times10^{-2}$ \\
\hline
\end{tabular}\\
\\
\textbf{Table 3: Comparison of numerical results  of example 6.1 with} $\textbf{h=0.025}$\\
\begin{tabular}{|c|c|c|c|c|c|c|}
\hline
\textbf{x} & \textbf{Our method} & \textbf{Runge-Kutta}  & \textbf{Predictor-} & \textbf{Multi-step} & \textbf{Day method}  \\
\multirow{2}{*}
 &$O(h^3)$  &  $O(h^3)$ & \textbf{Corrector} $O(h^3)$ & $O(h^3)$&   \\
\hline
  0.1 & $2.1\times10^{-6}$ & $ 3.3\times10^{-5}$ & $ 2.2\times10^{-6}$ & $7.7\times10^{-6}$& $3.3\times10^{-5}$ \\
\hline
  0.2 & $9.0 \times10^{-6}$ & $ 7.1\times10^{-5}$ & $ 1.1\times10^{-5}$& $2.6\times10^{-5}$& $7.1\times10^{-5}$\\
\hline
  0.3 & $2.1\times10^{-5}$ & $ 1.2\times10^{-4}$  & $2.7\times10^{-5}$ & $5.7\times10^{-5}$ & $1.1\times10^{-4}$\\
\hline
  0.4 & $4.1\times10^{-5}$ & $1.7\times10^{-4}$  & $5.0\times10^{-5}$ & $1.0\times10^{-4}$ & $1.7\times10^{-4}$\\
\hline
  0.5 & $7.0\times10^{-5}$ & $2.4\times10^{-4}$ & $8.4\times10^{-5}$ & $1.7\times10^{-4}$ & $2.5\times10^{-4}$ \\
\hline
 0.6 & $1.7\times10^{-4}$ & $ 3.3\times10^{-4}$  & $ 1.3\times10^{-4}$ & $2.6\times10^{-4}$ & $3.5\times10^{-4}$ \\
\hline
  0.7 & $1.8\times10^{-4}$ & $4.5 \times10^{-4} $ & $2.0\times10^{-4}$ & $4.0\times10^{-4}$ & $5.0\times10^{-4}$ \\
\hline
  0.8 & $2.8\times10^{-4}$ & $6.3 \times10^{-4} $ & $3.1\times10^{-4}$ & $5.9\times10^{-4}$ & $7.0\times10^{-4}$ \\
\hline
  0.9 & $4.3\times10^{-4}$ & $8.7 \times10^{-4} $ & $4.6\times10^{-4}$ & $8.8\times10^{-4}$ & $1.0\times10^{-3}$ \\
\hline
  1 & $ 6.6\times10^{-4}$ &  $1.2\times10^{-3} $ & $7.0\times10^{-4}$ & $1.3\times10^{-3}$ & $1.4\times10^{-3}$ \\
\hline
\end{tabular}

\begin{Ex} \cite{Dehghan}\label{ex2}
Consider the nonlinear Volterra integro-differential equation
\begin{equation}
y'(x) = 2x -\frac{1}{2} sin(x^4) + \int_{0}^x x^2 t cos(x^2y(t))dt,\quad y(0) = 0.\label{20}
\end{equation}
\end{Ex}
The exact solution of Eq.(\ref{20}) is $y(x) = x^2$. We compare the maximum absolute errors obtained using meshless method given in \cite{Dehghan} and CPU times used for different values of n in the interval $[0,1]$. It is observed that the error in our method is very less as compared with linear meshless method. Further, our method is more time efficient than both linear and quadratic meshless method.
\\
\\
\\
\\
\textbf{Table 4: Maximum absolute errors and CPU times used for different values of n.}\\
\begin{tabular}{|c|c|c|c|c|c|c|c|}
\hline
\textbf{n} & \textbf{Our method }& \textbf{Time}  &\textbf{ Linear Meshless}  & \textbf{Time }& \textbf{Quadratic Meshless} & \textbf{Time}\\
\multirow{2}{*}
 &  &  & \textbf{ method} &  & \textbf{method} &  \\
\hline
 5 & $1.49\times10^{-3}$& 0.00 & $ 2.83\times10^{-3}$ & 0.42 &$2.21\times10^{-4}$ & 0.58  \\
  \hline
    9 & $4.54\times10^{-4}$& 0.00 & $ 7.22\times10^{-4}$& 0.68 & $ 3.41\times10^{-5}$& 0.79 \\
  \hline
   17 & $1.26\times10^{-4}$ & 0.00 &$ 2.25\times10^{-4}$ & 0.76 & $9.98\times10^{-6}$ & 1.16 \\
  \hline
   33 & $3.35\times10^{-5}$ & 0.01 &$6.93\times10^{-5}$  & 0.87 &$3.26\times10^{-6}$& 1.85 \\
  \hline
  65 & $8.26\times10^{-6}$& 0.04 & $2.19\times10^{-5}$& 0.95 & $1.87\times10^{-6}$ & 2.15\\
  \hline
  129 & $2.19\times10^{-6}$& 0.15 & $ 6.58\times10^{-6}$& 1.40 & $1.28\times10^{-6}$& 2.91\\
  \hline   
\end{tabular}

\begin{Ex} \label{ex3}
Consider the VIDE \cite{Dehghan}
\begin{equation}
y'(x) = 1 -\frac{x}{2}+ \frac{x e^{-x^2}}{2} + \int_{0}^x x t e^{-y^2(t)} dt,\quad y(0) = 0.\label{21}
\end{equation}
\end{Ex}
The exact solution of Eq.(\ref{21}) is $y(x) = x$. We compare the maximum absolute errors obtained using meshless method given in \cite{Dehghan} and CPU time used for different values of n in the interval $[0,1]$. It is observed that our method is more accurate as compared with linear meshless and more time efficient than both linear and quadratic meshless method.\\
\textbf{Table 5: Maximum absolute errors and CPU times used for different values of n.}\\
\begin{tabular}{|c|c|c|c|c|c|c|c|}
\hline
\textbf{n} & \textbf{Our method }& \textbf{Time}  &\textbf{ Linear Meshless}  & \textbf{Time }& \textbf{Quadratic Meshless} & \textbf{Time}\\
\multirow{2}{*}
 &  &  & \textbf{ method} &  & \textbf{method} &  \\
\hline
 5 & $8.96\times10^{-3}$& 0.00 & $ 5.84\times10^{-3}$ & 0.45 &$ 1.96\times10^{-4}$ & 0.46  \\
  \hline
    9 & $2.33\times10^{-4}$& 0.00 & $ 1.75\times10^{-3}$& 0.57 & $ 2.42\times10^{-5}$& 0.58 \\
  \hline
   17 & $6.13\times10^{-5}$ & 0.00 &$ 4.88\times10^{-4}$ & 0.65 & $3.48\times10^{-6}$ & 0.73 \\
  \hline
   33 & $1.59\times10^{-5}$ & 0.01 &$1.30\times10^{-4}$  & 0.72 &$6.55\times10^{-7}$& 0.87 \\
  \hline
  65 & $4.10\times10^{-6}$& 0.04 & $3.21\times10^{-5}$& 1.68 & $5.75\times10^{-7}$ & 1.66\\
  \hline
  129 & $1.04\times10^{-6}$& 0.17 & $ 8.19\times10^{-6}$& 4.83 & $ 6.05\times10^{-7}$& 5.61\\
  \hline
\end{tabular}
\\
\\
We compare our solution with exact solution for $h=0.01$  in Fig.3. It is observed that our solution coincides with exact solution.

\begin{center}
\begin{tabular}{c}
\includegraphics[scale=1.0]{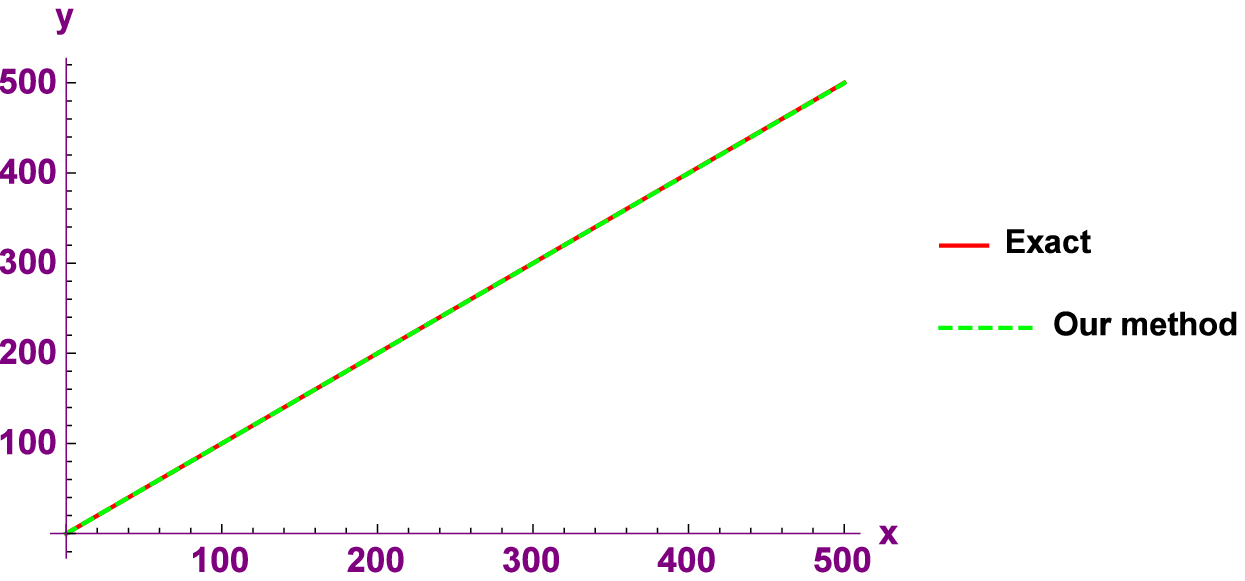}\\
Fig.3: Comparison of solutions of Eq.(\ref{21}) for $h=0.01$.\\
\end{tabular}
\end{center}

\begin{Ex} \label{ex4}
Consider the  VIDE 
\begin{equation}
y'(x) = 1 + \int_{0}^x e^{-t} y^2(t) dt,\quad y(0) = 1.\label{23}
\end{equation}
\end{Ex}
The exact solution of Eq.(\ref{23}) is $y(x) = e^x$. We compare our solution with exact solution for $h=0.01$  in Fig.4. It is observed that our solution is well in agreement  with exact solution.
\begin{center}
\begin{tabular}{c}
\includegraphics[scale=1.0]{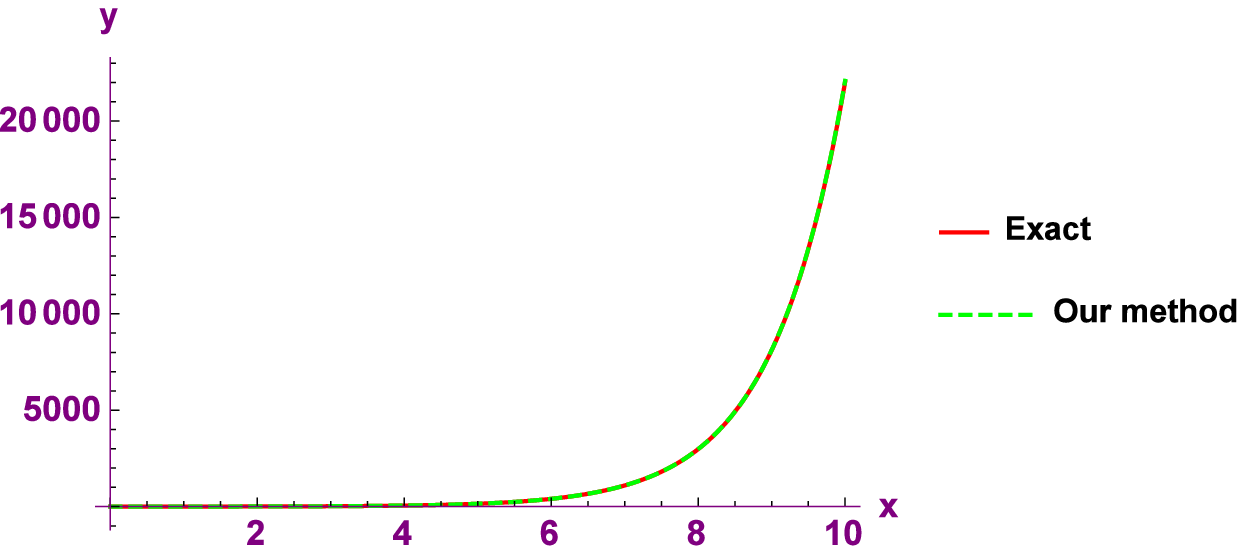}\\
Fig.4: Comparison of solutions of Eq.(\ref{23}) for $h=0.01$.\\
\end{tabular}
\end{center}

\section{Software Package}\label{Package}
In this section we provide software package based on Mathematica-$10$ to solve VIDEs using numerical method.
\par This software is used for solving VIDEs of the form (\ref{21}). One has to provide the value of $f(x,y)$ and $K(x,t,y)$ in first two windows, initial condition $y_0$, step size $h$ and number of steps $n$ in third, fourth and fifth windows respectively. The last window gives required solution curve. We solve Eq.(\ref{20}) using this software package as shown in Fig.5.  \\
\begin{center}
\begin{tabular}{c}
\includegraphics[scale=1]{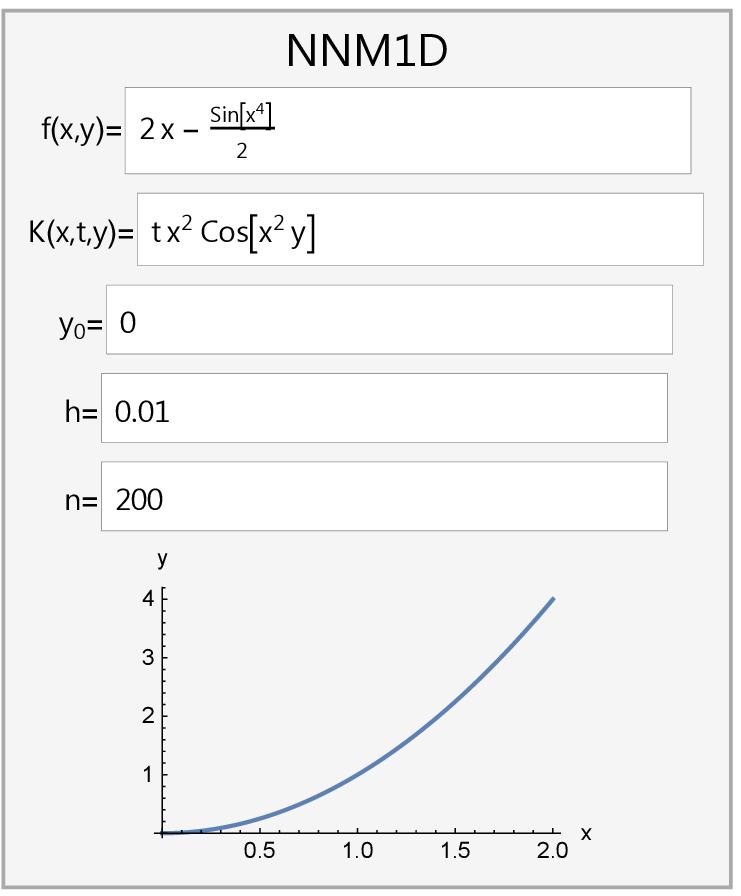}\\
Fig.5: Software Package\\
\end{tabular}
\end{center}

\section{Conclusions}\label{concl}
A new method proposed in this paper is a combination of trapezium rule and DJM. The DJM is used to find the approximation to the unknown term on the right side of the difference formula obtained by using trapezium rule. This combination is proved to be very efficient. The error in this third order method is very less as compared to other similar methods. However, the time taken for the computation was very less. The stability and bifurcation analysis also show the power of this method. An ample number of examples show the wide range of applicability of this method.
\\ 

\textbf{Acknowledgements:}\\
S. Bhalekar acknowledges CSIR, New Delhi for funding through Research Project [25(0245)/15/EMR-II].

\end{document}